\documentclass{elsart}
\journal{J. Comput. Appl. Math.}
\usepackage{amssymb}
\usepackage{amsmath}
\newcommand{\D}{\mathcal{D}}
\newtheorem{theorem}{Theorem}[section]
\numberwithin{equation}{section}
\begin{document}
\begin{frontmatter}
\title{Multiple little $q$-Jacobi polynomials\thanksref{label1}}
\thanks[label1]{This work was supported by INTAS Research Network 03-51-6631 and FWO projects G.0184.02 and G.0455.04}
\author{Kelly Postelmans and Walter Van Assche}
\ead{\{Kelly.Postelmans,Walter.VanAssche\}@wis.kuleuven.ac.be}
\address{Katholieke Universiteit Leuven, Department of Mathematics,\\
Celestijnenlaan 200B, B-3001 Leuven, Belgium}

\begin{abstract}
We introduce two kinds of multiple little $q$-Jacobi polynomials $p_{\vec{n}}$ with
multi-index $\vec{n}=(n_1,n_2,\ldots,n_r)$ and degree $|\vec{n}| = n_1+n_2+\cdots+n_r$
by imposing orthogonality conditions with respect to $r$ discrete little $q$-Jacobi measures
on the exponential lattice $\{q^k, k=0,1,2,3,\ldots\}$, where $0 < q <1$.
We show that
these multiple little $q$-Jacobi polynomials have useful $q$-difference properties,
such as a Rodrigues formula (consisting of a product of $r$ difference operators). Some properties of the zeros of these polynomials and some
asymptotic properties will be given as well.
\end{abstract}

\end{frontmatter}

\section{Little $q$-Jacobi polynomials}
Little $q$-Jacobi polynomials are orthogonal polynomials on the
exponential lattice $\{q^k,\ k=0,1,2,\ldots\}$, where $0 < q < 1$.
In order to express the orthogonality relations, we will use the
$q$-integral
\begin{equation} \label{eq:qint}
   \int_0^1 f(x)\, d_qx = (1-q) \sum_{k=0}^\infty q^k f(q^k),
\end{equation}
(see, e.g., \cite[\S 10.1]{AAR}, \cite[\S 1.11]{GR}) where $f$ is
a function on $[0,1]$ which is continuous at $0$. The
orthogonality is given by
\begin{equation} \label{eq:ortho}
  \int_0^1 p_n(x;\alpha,\beta|q) x^k w(x;\alpha,\beta|q)\, d_qx = 0, \qquad k=0,1,\ldots,n-1,
\end{equation}
where
\begin{equation} \label{eq:w}
   w(x;a,b|q) = \frac{(qx;q)_\infty}{(q^{\beta+1}x;q)_\infty} x^\alpha .
\end{equation}
We have used the notation
\[  (a;q)_n = \prod_{k=0}^{n-1} (1-aq^k), \qquad (a;q)_\infty = \prod_{k=0}^\infty (1-aq^k). \]
In order that the $q$-integral of $w$ is finite, we need to impose the restrictions $\alpha,\beta > -1$. The orthogonality conditions (\ref{eq:ortho}) determine the polynomials $p_n(x;\alpha,\beta|q)$ up to
a multiplicative factor. In this paper we will always use monic polynomials and these are uniquely determined by the orthogonality conditions.
The $q$-binomial theorem
\begin{equation} \label{eq:qbin}
   \sum_{n=0}^{\infty} \frac{(a;q)_n}{(q;q)_n} z^n = \frac{(az;q)_\infty}{(z;q)_\infty}, \qquad |z|, |q| < 1,
\end{equation}
(see, e.g., \cite[\S 10.2]{AAR}, \cite[\S 1.3]{GR}) implies that
\[  \lim_{q \to 1} w(x;\alpha,\beta|q) = (1-x)^\beta x^\alpha, \qquad 0 < x < 1,
\]
so that $w(x;\alpha,\beta|q)$ is a $q$-analog of the beta density on
$[0,1]$, and hence
\[  \lim_{q \to 1} p_n(x;\alpha,\beta|q) = P_n^{(\alpha,\beta)}(x) , \]
where $P_n^{(\alpha,\beta)}$ are the monic Jacobi polynomials on $[0,1]$.
Little $q$-Jacobi polynomials appear in representations of quantum $SU(2)$
\cite{K}, \cite{MMNNU}, and the special case of little $q$-Legendre polynomials was used to prove irrationality of a $q$-analog of the
harmonic series and $\log 2$ \cite{V}. Their role in partitions was
described in \cite{AA}. A detailed list of formulas for the little
$q$-Jacobi polynomials can be found in \cite[\S 3.12]{KS}, but note that
in that reference the polynomial $p_n(x;a,b|q)$ is not monic and that $a=q^\alpha, b=q^\beta$.
Useful formulas are the \textit{lowering operation}
\begin{equation} \label{eq:low}
    \D_q p_n(x;\alpha,\beta|q) = \frac{1-q^n}{1-q} p_{n-1}(x;\alpha+1,\beta+1|q),
\end{equation}
where $\D_q$ is the $q$-difference operator
\begin{equation} \label{eq:qdif}
   \D_qf(x) = \begin{cases}
            \displaystyle \frac{f(x)-f(qx)}{(1-q)x}, & \textrm{if $x \neq 0$}, \\
            f'(0), & \textrm{if $x=0$},
        \end{cases}
\end{equation}
and the \textit{raising operation}
\begin{multline}  \label{eq:rais}
   \D_p [w(x;\alpha,\beta|q) p_n(x;\alpha,\beta|q)] \\
   = - \frac{1-q^{n+\alpha+\beta}}{(1-q)q^{n+\alpha-1}} w(x;\alpha-1,\beta-1|q)
  p_{n+1}(x;\alpha-1,\beta-1|q),
\end{multline}
where $p=1/q$.
Repeated application of the raising operator gives the \textit{Rodrigues formula}
\begin{equation}  \label{eq:Rod}
  w(x;\alpha,\beta|q)p_n(x;\alpha,\beta|q) = \frac{(-1)^n(1-q)^nq^{\alpha n+n(n-1)}}{(q^{\alpha+\beta+n+1};q)_n} \D_p^n w(x;\alpha+n,\beta+n|q).
\end{equation}
A combination of the raising and the lowering operation gives a \textit{second order
$q$-difference equation}.
The Rodrigues formula enables us to give an explicit expression as a basic hypergeometric
sum:
\[   p_n(x;\alpha,\beta|q) = \frac{x^n q^{n(n+\alpha)} (q^{-n-\alpha};q)_n}{(q^{n+\alpha+\beta+1};q)_n}  {}_3\phi_2 \left( \left. \begin{array}{c}
   q^{-n},q^{-n-\alpha},1/x \\ q^{\beta+1}, 0
   \end{array} \right| q;q \right),  \]
which by some elementary transformations can also be written as
\begin{eqnarray}
  p_n(x;\alpha,\beta|q) &=& \frac{q^{(n+\alpha)n}(q^{-n-\alpha};q)_n}{(q^{n+\alpha+\beta+1};q)_n}
   {}_2\phi_1 \left( \left. \begin{array}{c}
   q^{-n},q^{n+\alpha+\beta+1} \\ q^{\alpha+1}
   \end{array} \right| q;qx \right) \nonumber \\
  & = & \frac{q^{(n+\alpha)n}(q^{-n-\alpha};q)_n}{(q^{n+\alpha+\beta+1};q)_n}
  \sum_{k=0}^n \frac{(q^{-n};q)_k (q^{n+\alpha+\beta+1};q)_k}{(q^{\alpha+1};q)_k} \frac{q^kx^k}{(q;q)_k} .  \label{eq:expl}
\end{eqnarray}

\section{Multiple orthogonal polynomials}
Multiple orthogonal polynomials (of type II) are polynomials satisfying orthogonality
conditions with respect to $r \geq 1$ positive measures \cite{A} \cite{ACV} \cite[\S 4.3]{NS} \cite{VC}.
Let $\mu_1, \mu_2, \ldots, \mu_r$
be $r$ positive measures on the real line and let $\vec{n} = (n_1,n_2,\ldots,n_r) \in \mathbb{N}^r$ be a multi-index of length $|\vec{n}|=n_1+n_2+\cdots+ n_r$.
The corresponding type II multiple orthogonal polynomial $p_{\vec{n}}$ is
a polynomial of degree $\leq |\vec{n}|$ satisfying the orthogonality relations
\begin{equation*}
  \int p_{\vec{n}}(x) x^k \, d\mu_j(x)  =  0, \qquad k=0,1,\ldots,n_j-1,\ j=1,2,\ldots,r.
\end{equation*}
These orthogonality relations give $|\vec{n}|$ homogeneous equations for the
$|\vec{n}|+1$ unknown coefficients of $p_{\vec{n}}$. We say that $\vec{n}$ is a normal
index if the orthogonality relations determine the polynomial $p_{\vec{n}}$ up to a multiplicative factor. Multiple orthogonal polynomials of type I (see, e.g., \cite{A} \cite[\S 4.3]{NS} \cite{ACV} \cite{VC}) will not be considered
in this paper. Multiple little $q$-Jacobi polynomials are
multiple orthogonal polynomials where the measures $\mu_1,\ldots,\mu_r$ are
supported on the exponential lattice $\{q^k,\ k=0,1,2,\ldots\}$ and are all of the form
$d\mu_i(x) = w(x;\alpha_i,\beta_i|q)\, d_qx$, where $w(x;\alpha,\beta|q)\, d_qx$ is the orthogonality measure
for little $q$-Jacobi polynomials. It turns out that in order to have formulas and identities similar
to those of the usual little $q$-Jacobi polynomials one needs to
keep one of the parameters $\alpha_i$ or $\beta_i$ fixed and  change the other parameters
for the $r$ measures. This gives two kinds of multiple little $q$-Jacobi polynomials.
Note that these multiple little $q$-Jacobi polynomials should not be confused
with multivariable little $q$-Jacobi polynomials, introduced by Stokman \cite{S}.
In \cite{PV} the multiple little $q$-Jacobi polynomials of the first kind are used
to prove some irrationality results for $\zeta_q(1)$ and $\zeta_q(2)$.

\subsection{Multiple little $q$-Jacobi polynomials of the first kind}
Multiple little $q$-Jacobi polynomials of the first kind $p_{\vec{n}}(x;\vec{\alpha},\beta|q)$
are monic polynomials of degree $|\vec{n}|$ satisfying the orthogonality relations
\begin{equation}  \label{eq:1kind}
  \int_0^1 p_{\vec{n}}(x;\vec{\alpha},\beta|q) x^k w(x;\alpha_j,\beta|q)\, d_qx  =  0, \qquad k=0,1,\ldots,n_j-1,\ j=1,2,\ldots,r,
\end{equation}
where $\alpha_1,\ldots,\alpha_r,\beta > -1$. Observe that all the measures are orthogonality
measures for little $q$-Jacobi polynomials with the same parameter $\beta$ but with different parameters $\alpha_j$. All the multi-indices will be normal when we impose the condition that $\alpha_i-\alpha_j \notin \mathbb{Z}$
whenever $i\neq j$, because then all the measures are absolutely continuous with respect to
$w(x;0,\beta|q)\, d_qx$ and the system of functions
\[   x^{\alpha_1},  x^{\alpha_1+1},\ldots,x^{\alpha_1+n_1-1},
      x^{\alpha_2},  x^{\alpha_2+1},\ldots,x^{\alpha_2+n_2-1}, \ldots,
      x^{\alpha_r},  x^{\alpha_r+1},\ldots,x^{\alpha_r+n_r-1}
\]
is a Chebyshev system on $(0,1)$, so that the measures $(\mu_1,\ldots,\mu_r)$ form a
so-called AT-system, which implies that all the multi-indices are normal \cite[Theorem 4.3]{NS}.

There are $r$ \textit{raising operations} for these multiple orthogonal polynomials.
\begin{theorem}
Suppose that $\alpha_1,\ldots,\alpha_r,\beta > 0$, with $\alpha_i-\alpha_j \notin \mathbb{Z}$
whenever $i \neq j$, and put $p=1/q$, then
\begin{multline}  \label{eq:rais1}
  \D_p \left[ w(x;\alpha_j,\beta|q) p_{\vec{n}}(x;\vec{\alpha},\beta|q) \right] \\
  =     \frac{q^{\alpha_j+\beta+|\vec{n}|}-1}{(1-q)q^{\alpha_j+|\vec{n}|-1}}
      w(x;\alpha_j-1,\beta-1|q)
   p_{\vec{n}+\vec{e}_j}(x;\vec{\alpha}-\vec{e}_j,\beta-1|q),
\end{multline}
for $1 \leq j \leq r$,
where $\vec{e}_1=(1,0,0,\ldots,0),\ldots,\vec{e}_r=(0,\ldots,0,0,1)$ are the standard unit vectors.
\end{theorem}
Observe that these operations raise one of the indices in the multi-index and lower the parameter $\beta$ and
one of the components of $\vec{\alpha}$.

\textbf{Proof:}
First observe that
\begin{multline*}
  \D_p\left[ w(x;\alpha_j,\beta|q)p_{\vec{n}}(x;\vec{\alpha},\beta|q) \right] \\
   = w(x;\alpha_j-1,\beta-1|q) \frac{(1-q^{\beta} x)p_{\vec{n}}(x;\vec{\alpha},\beta|q)-
   p^{\alpha_j}(1-x) p_{\vec{n}}(px;\vec{\alpha},\beta|q)}{1-p},
\end{multline*}
so that
\begin{equation} \label{eq:Q1}
 \D_p \left[ w(x;\alpha_j,\beta|q)p_{\vec{n}}(x;\vec{\alpha},\beta|q) \right]
   = - \frac{1-q^{\alpha_j+\beta+|\vec{n}|}}{(1-q)q^{\alpha_j+|\vec{n}|-1}} w(x;\alpha_j-1,\beta-1|q) Q_{|\vec{n}|+1}(x),
\end{equation}
where $Q_{|\vec{n}|+1}$ is a monic polynomial of degree $|\vec{n}|+1$. We will show that this monic polynomial
$Q_{|\vec{n}|+1}$ satisfies the multiple orthogonality conditions (\ref{eq:1kind}) of
 $p_{\vec{n}+\vec{e}_j}(x;\vec{\alpha}-\vec{e}_j,\beta-1|q)$ and hence, since all $\alpha_i-\alpha_j \notin \mathbb{Z}$ whenever $i\neq j$, the unicity of the
multiple orthogonal polynomials implies that
$Q_{|\vec{n}|+1}(x) = p_{\vec{n}+\vec{e}_j}(x;\vec{\alpha}-\vec{e}_j,\beta-1|q)$.
Integration by parts for the $q$-integral is given by the rule
\begin{equation}  \label{eq:parts}
  \int_0^1 f(x) \D_p g(x) = -q \int_0^1 g(x)\D_q f(x), \qquad \textrm{if } g(p)=0.
\end{equation}
If we apply this, then
\begin{multline*}
  \frac{1-q^{\alpha_j+\beta+|\vec{n}|}}{(1-q)q^{\alpha_j+|\vec{n}|-1}}
  \int_0^1 x^k w(x;\alpha_j-1,\beta-1|q) Q_{|\vec{n}|+1}(x)\, d_qx   \\
 =   - q \int_0^1 w(x;\alpha_j,\beta|q) p_{\vec{n}}(x;\vec{\alpha},\beta|q) \D_q x^k \, d_qx ,
\end{multline*}
and since
\[   \D_q x^k = \begin{cases}   \frac{1-q^k}{1-q} x^{k-1} & \textrm{if $k \geq 1$}, \\
                    0 & \textrm{if $k=0$,}
               \end{cases}  \]
we find that
\[  \int_0^1 x^k w(x;\alpha_j-1,\beta-1|q) Q_{|\vec{n}|+1}(x)\, d_qx = 0, \qquad k=0,1,\ldots,n_j. \]
For the other components $\alpha_i$ $(i \neq j)$ of $\vec{\alpha}$ we have
\begin{eqnarray*}
  \lefteqn{ \frac{1-q^{\alpha_j+\beta+|\vec{n}|}}{(1-q)q^{\alpha_j+|\vec{n}|-1}}
  \int_0^1 x^k w(x;\alpha_i,\beta-1|q) Q_{|\vec{n}|+1}(x)\, d_qx} \quad & & \\
& = &    \frac{1-q^{\alpha_j+\beta+|\vec{n}|}}{(1-q)q^{\alpha_j+|\vec{n}|-1}}
  \int_0^1 x^{k+\alpha_i-\alpha_j+1} w(x;\alpha_j-1,\beta-1|q) Q_{|\vec{n}|+1}(x)\, d_qx  \\
& = &   -\ q \int_0^1 w(x;\alpha_j,\beta|q) p_{\vec{n}}(x;\vec{\alpha},\beta|q) \D_q x^{k+\alpha_i-\alpha_j+1} \, d_qx ,
\end{eqnarray*}
and since $\alpha_i-\alpha_j \notin \mathbb{Z}$ we have
\[   \D_q x^{k+\alpha_i-\alpha_j+1} = \frac{1-q^{k+\alpha_i-\alpha_j}}{1-q} x^{k+\alpha_i-\alpha_j}, \]
hence
\[  \int_0^1 x^k w(x;\alpha_i,\beta-1|q) Q_{|\vec{n}|+1}(x)\, d_qx = 0, \qquad k=0,1,\ldots,n_i-1. \]
Hence all the orthogonality conditions for $p_{\vec{n}+\vec{e}_j}(x;\vec{\alpha}-\vec{e}_j,\beta-1|q)$ are
indeed satisfied. \qed

As a consequence we find a \textit{Rodrigues formula}:

\begin{theorem}
The multiple little $q$-Jacobi polynomials of the first kind are given by
\begin{equation}  \label{eq:RodI}
 p_{\vec{n}}(x;\vec{\alpha},\beta|q)  = C(\vec{n},\vec{\alpha},\beta)
 \frac{(q^{\beta+1}x;q)_\infty}{(qx;q)_\infty}
   \prod_{j=1}^r \left( x^{-\alpha_j} \D_p^{n_j} x^{\alpha_j+n_j} \right)
   \frac{(qx;q)_\infty}{(q^{\beta+|\vec{n}|+1}x;q)_\infty},
\end{equation}
where the product of the difference operators can be taken in any order and
\[   C(\vec{n},\vec{\alpha},\beta) =  (-1)^{|\vec{n}|}
      \frac{(1-q)^{|\vec{n}|} q^{\sum_{j=1}^r (\alpha_j-1) n_j + \sum_{1\leq j\leq k \leq r}
 n_jn_k}}{ \prod_{j=1}^r (q^{\alpha_j+\beta+|\vec{n}|+1};q)_{n_j}}. \]
\end{theorem}

\textbf{Proof:}
If we apply the  raising operator for $\alpha_j$ recursively $n_j$ times, then
\begin{multline}  \label{eq:raisn}
 \D_p^{n_j} w(x;\alpha_j,\beta|q)p_{\vec{m}}(x;\vec{\alpha},\beta|q)
  = (-1)^{n_j} \frac{(q^{\alpha_j+\beta+|\vec{m}|-n_j+1};q)_{n_j}}
  {(1-q)^{n_j}q^{(\alpha_j+|\vec{m}|-1)n_j}} \\
\times  w(x;\alpha_j-n_j,\beta-n_j|q)
p_{\vec{m}+n_j\vec{e}_j}(x;\vec{\alpha}-n_j\vec{e}_j,\beta-n_j|q).
\end{multline}
Use this expression with $\vec{m}=\vec{0}$ and $j=1$, then
\begin{multline*}
    \D_p^{n_1} w(x;\alpha_1,\beta|q) = (-1)^{n_1}
  \frac{(q^{\alpha_1+\beta-n_1+1};q)_{n_1}}{(1-q)^{n_1} q^{(\alpha_1-1)n_1}} \\
 \times w(x;\alpha_1-n_1,\beta-n_1|q)
p_{n_1\vec{e}_1}(x;\vec{\alpha}-n_1\vec{e}_1,\beta-n_1|q).
\end{multline*}
Multiply both sides by $w(x;\alpha_2,\beta-n_1|q)$ and divide by
$w(x;\alpha_1-n_1,\beta-n_1|q)$, then
\begin{multline*}
   x^{n_1+\alpha_2-\alpha_1} \D_p^{n_1} w(x;\alpha_1,\beta|q) = (-1)^{n_1}
 \frac{(q^{\alpha_1+\beta-n_1+1};q)_{n_1}}{(1-q)^{n_1} q^{(\alpha_1-1)n_1}} \\
\times  w(x;\alpha_2,\beta-n_1|q)
p_{n_1\vec{e}_1}(x;\vec{\alpha}-n_1\vec{e}_1,\beta-n_1|q).
\end{multline*}
Apply (\ref{eq:raisn}) with $j=2$, then
\begin{multline*}
  \D_p^{n_2} x^{n_1+\alpha_2-\alpha_1} \D_p^{n_1} w(x;\alpha_1,\beta|q) = (-1)^{n_1+n_2}
  \frac{(q^{\alpha_1+\beta-n_1+1};q)_{n_1}(q^{\alpha_2+\beta-n_2+1};q)_{n_2}}
 {(1-q)^{n_1+n_2} q^{(\alpha_1-1)n_1+(\alpha_2-1+n_1)n_2}} \\
 \times   w(x;\alpha_2-n_2,\beta-n_1-n_2|q)
p_{n_1\vec{e}_1+n_2\vec{e}_2}(x;\vec{\alpha}-n_1\vec{e}_1-n_2\vec{e}_2,\beta-n_1-n_2|q).
\end{multline*}
Continuing this way we arrive at
\begin{multline*}
\left(\D_p^{n_r} x^{\alpha_r}\right) \left( x^{n_{r-1}-\alpha_{r-1}} \D_p^{n_{r-1}} x^{\alpha_{r-1}} \right)
  \cdots \left( x^{n_1-\alpha_1} \D_p^{n_1} \right) w(x;\alpha_1,\beta|q)\\
 =  \frac{(-1)^{|\vec{n}|}\prod_{j=1}^r (q^{\alpha_j+\beta-n_j+1};q)_{n_j}}
{(1-q)^{|\vec{n}|} q^{\sum_{j=1}^r (\alpha_j-1)n_j + \sum_{1 \leq j < k \leq r} n_jn_k}}
   w(x;\alpha_r-n_r,\beta-|\vec{n}||q)
  p_{\vec{n}}(x;\vec{\alpha}-\vec{n},\beta-|\vec{n}||q).
\end{multline*}
Now replace each $\alpha_j$ by $\alpha_j+n_j$ and $\beta$ by $\beta+|\vec{n}|$, then
the required expression follows.
 The order in which we took the raising operators
is irrelevant.
\qed

 We can obtain an explicit expression of the multiple little $q$-Jacobi polynomials of
the first kind using this Rodrigues formula. Indeed, if we use the $q$-binomial theorem,
then
\[   \frac{(qx;q)_\infty}{(q^{\beta+|\vec{n}|+1}x;q)_\infty}
   = \sum_{k=0}^\infty \frac{(q^{-\beta-|\vec{n}|};q)_k}{(q;q)_k}
  q^{(\beta+|\vec{n}|+1)k} x^k. \]
Use this in (\ref{eq:RodI}), together with
\[  x^{-\alpha} \D_p^n x^{\alpha+n+k} =
  \frac{(q^{\alpha+1};q)_n}{(1-q)^n} \frac{(q^{\alpha+n+1};q)_k}{(q^{\alpha+1};q)_k}
    q^{-n(k+\alpha)-n(n-1)/2} x^k, \]
then this gives
\begin{multline}  \label{eq:infinsumI}
 p_{\vec{n}}(x;\vec{\alpha},\beta|q) = C(\vec{n},\vec{\alpha},\beta)
 \frac{\prod_{j=1}^r (q^{\alpha_j+1};q)_{n_j}}{(1-q)^{|\vec{n}|}} q^{-\sum_{j=1}^r \alpha_jn_j
  -\sum_{j=1}^r \binom{n_j}{2} }\\
 \frac{(q^{\beta+1}x;q)_\infty}{(qx;q)_\infty}
  {}_{r+1}\phi_r \left( \left. \begin{array}{c}
     q^{-\beta-|\vec{n}|}, q^{\alpha_1+n_1+1}, \ldots, q^{\alpha_r+n_r+1} \\
        q^{\alpha_1+1},\ldots,q^{\alpha_r+1} \end{array} \right| q;q^{\beta+1}x
   \right).
\end{multline}
This explicit expression uses a non-terminating basic hypergeometric series, except when $\beta$ is an
integer. Another representation, using only finite sums, can be obtained by using the Rodrigues
formula (\ref{eq:Rod}) $r$ times. For $r=2$ this gives
\begin{theorem}
The multiple little $q$-Jacobi polynomials of the first kind (for $r=2$) are given by
\begin{multline} \label{eq:explI}
  p_{n,m}(x;(\alpha_1,\alpha_2),\beta|q) = \frac{q^{nm +m^2+n^2+\alpha_1n+\alpha_2m}
   (q^{-\alpha_1-n};q)_n(q^{-\alpha_2-m};q)_m}{(q^{\alpha_1+\beta+n+m+1};q)_n(q^{\alpha_2+\beta+n+m+1};q)_m} \\
\times  \sum_{\ell=0}^n \sum_{k=0}^m \frac{(q^{-n};q)_\ell (q^{-m};q)_k (q^{\alpha_2+\beta+m+n+1};q)_k
(q^{\alpha_1+\beta+n+1};q)_{k+\ell} (q^{\alpha_1+n+1};q)_k}{(q^{\alpha_2+1};q)_k(q^{\alpha_1+1};q)_{k+\ell}
(q^{\alpha_1+\beta+n+1};q)_k} \\
 \times \frac{q^{k+\ell} x^{k+\ell}}{q^{kn}(q;q)_k (q;q)_\ell} .
\end{multline}
\end{theorem}

\textbf{Proof:}
For $r=2$ the Rodrigues formula (\ref{eq:RodI}) is
\begin{multline*}
  p_{n,m}(x;(\alpha_1,\alpha_2),\beta|q) = \frac{(-1)^{n+m}(1-q)^{n+m}q^{\alpha_1n+\alpha_2m-n-m+nm+n^2+m^2}}
   {(q^{\alpha_1+\beta+n+m+1};q)_n (q^{\alpha_2+\beta+n+m+1};q)_m} \\
\times \frac{(q^{\beta+1}x;q)_\infty}{(qx;q)_\infty} x^{-\alpha_1} \D_p^n x^{\alpha_1+n-\alpha_2} \D_p^m
   x^{\alpha_2+m} \frac{(qx;q)_\infty}{(q^{\beta+n+m+1}x;q)_\infty} .
\end{multline*}
Observe that by the Rodrigues formula (\ref{eq:Rod}) for the little $q$-Jacobi polynomials
\begin{multline*}   \D_p^m x^{\alpha_2+m} \frac{(qx;q)_\infty}{(q^{\beta+n+m+1}x;q)_\infty} \\
  = \frac{(-1)^m (q^{\alpha_2+\beta+n+m+1};q)_m}{(1-q)^mq^{\alpha_2m+m^2-m}} \frac{(qx;q)_\infty}{(q^{\beta+n+1}x;q)_\infty}
   x^{\alpha_2} p_m(x;\alpha_2,\beta+n|q),
\end{multline*}
and hence
\begin{multline*}
  p_{n,m}(x;(\alpha_1,\alpha_2),\beta|q) = \frac{(-1)^{n}(1-q)^{n}q^{\alpha_1n-n+nm+n^2}}
   {(q^{\alpha_1+\beta+n+m+1};q)_n} \\
\times \frac{(q^{\beta+1}x;q)_\infty}{(qx;q)_\infty} x^{-\alpha_1} \D_p^n x^{\alpha_1+n}
   \frac{(qx;q)_\infty}{(q^{\beta+n+1}x;q)_\infty} p_m(x;\alpha_2,\beta+n|q) .
\end{multline*}
Now use the explicit expression (\ref{eq:expl}) to find
\begin{multline*}
  p_{n,m}(x;(\alpha_1,\alpha_2),\beta|q) = \frac{(-1)^{n}(1-q)^{n}q^{\alpha_1n+\alpha_2m-n+nm+n^2+m^2}
  (q^{-m-\alpha_2};q)_m}{(q^{\alpha_1+\beta+n+m+1};q)_n (q^{\alpha_2+\beta+n+m+1};q)_m} \\
\times \frac{(q^{\beta+1}x;q)_\infty}{(qx;q)_\infty} x^{-\alpha_1}
 \sum_{k=0}^m \frac{(q^{-m};q)_k(q^{\alpha_2+\beta+n+m+1};q)_k q^k}{(q^{\alpha_2+1};q)_k(q;q)_k}
\D_p^n x^{\alpha_1+n+k}    \frac{(qx;q)_\infty}{(q^{\beta+n+1}x;q)_\infty} .
\end{multline*}
In this expression we recognize
\begin{multline*}  \D_p^n x^{\alpha_1+n+k} \frac{(qx;q)_\infty}{(q^{\beta+n+1}x;q)_\infty} \\
    = \frac{(-1)^n(q^{\alpha_1+\beta+k+n+1};q)_n}{(1-q)^n q^{\alpha_1n+kn+n^2-n}} x^{\alpha_1+k}
   \frac{(qx;q)_\infty}{(q^{\beta+1}x;q)_\infty} p_n(x;\alpha_1+k,\beta|q),
\end{multline*}
hence
\begin{multline*}
  p_{n,m}(x;(\alpha_1,\alpha_2),\beta|q) = \frac{q^{\alpha_2m+nm+m^2}
  (q^{-m-\alpha_2};q)_m}{(q^{\alpha_1+\beta+n+m+1};q)_n (q^{\alpha_2+\beta+n+m+1};q)_m} \\
\times  \sum_{k=0}^m \frac{(q^{-m};q)_k(q^{\alpha_2+\beta+n+m+1};q)_k (q^{\alpha_1+\beta+k+n+1};q)_n q^k}
  {(q^{\alpha_2+1};q)_k(q;q)_kq^{kn}}
    x^k p_n(x;\alpha_1+k,\beta|q).
\end{multline*}
If we use the explicit expression (\ref{eq:expl}) for the little $q$-Jacobi polynomials once more, then
after some simplifications we finally arrive at (\ref{eq:explI}). \qed

\subsection{Multiple little $q$-Jacobi polynomials of the second kind}
Multiple little $q$-Jacobi polynomials of the second kind $p_{\vec{n}}(x;\alpha,\vec{\beta}|q)$
are monic polynomials of degree $|\vec{n}|$ satisfying the orthogonality relations
\begin{equation} \label{eq:2kind}
  \int_0^1 p_{\vec{n}}(x;\alpha,\vec{\beta}|q) x^k w(x;\alpha,\beta_j|q)\, d_qx  =  0,
  \qquad k=0,1,\ldots,n_j-1,\ j=1,2,\ldots,r,
\end{equation}
where $\alpha,\beta_1,\ldots,\beta_r > -1$. Observe that all the measures are orthogonality
measures for little $q$-Jacobi polynomials with the same parameter $\alpha$ but with different parameters $\beta_j$. All the multi-indices will be normal when we impose the condition that $\beta_i-\beta_j \notin \mathbb{Z}$
whenever $i\neq j$, because then all the measures are absolutely continuous with respect to
$(qx;q)_\infty w(x;\alpha,0|q)\, d_qx$ and the system of functions
\begin{multline*}
   \frac{1}{(q^{\beta_1+1}x;q)_\infty},  \frac{x}{(q^{\beta_1+1}x;q)_\infty} ,\ldots,\frac{x^{n_1-1}}{(q^{\beta_1+1}x;q)_\infty},
 \frac{1}{(q^{\beta_2+1}x;q)_\infty},  \frac{x}{(q^{\beta_2+1}x;q)_\infty} ,  \\
  \ldots,\frac{x^{n_2-1}}{(q^{\beta_2+1}x;q)_\infty}, \ldots,
  \frac{1}{(q^{\beta_r+1}x;q)_\infty},  \frac{x}{(q^{\beta_r+1}x;q)_\infty} ,\ldots,\frac{x^{n_r-1}}{(q^{\beta_r+1}x;q)_\infty}
\end{multline*}
is a Chebyshev system\footnote{The fact that this system is a Chebyshev system
is not obvious but is left as an advanced problem for the reader.}
 on $[0,1]$,
so that the vector of measures $(\mu_1,\ldots,\mu_r)$ forms an
AT-system, which implies that all the  multi-indices are normal
\cite[Theorem 4.3]{NS}.

Again there are $r$ \textit{raising operations}
\begin{theorem}
Suppose that $\alpha,\beta_1,\ldots,\beta_r > 0$, with $\beta_i-\beta_j \notin \mathbb{Z}$
when $i \neq j$, and put $p=1/q$, then
\begin{multline}  \label{eq:rais2}
  \D_p \left[ w(x;\alpha,\beta_j|q) p_{\vec{n}}(x;\alpha,\vec{\beta}|q) \right] \\
  =     \frac{q^{\alpha+\beta_j+|\vec{n}|}-1}{(1-q)q^{\alpha+|\vec{n}|-1}}
     w(x;\alpha-1,\beta_j-1|q)
   p_{\vec{n}+\vec{e}_j}(x;\alpha-1,\vec{\beta}-\vec{e}_j|q),
\end{multline}
for $1 \leq j \leq r$,
where $\vec{e}_1=(1,0,0,\ldots,0),\ldots,\vec{e}_r=(0,\ldots,0,0,1)$ are the standard unit vectors.
\end{theorem}
Observe that these operations raise one of the indices in the multi-index and lower the parameter $\alpha$ and
one of the components of $\vec{\beta}$.

\textbf{Proof:}
Again we see that
\begin{equation} \label{eq:Q2}
 \D_p \left[ w(x;\alpha,\beta_j|q)p_{\vec{n}}(x;\alpha,\vec{\beta}|q) \right]
   = \frac{q^{\alpha+\beta_j+|\vec{n}|}-1}{(1-q)q^{\alpha+|\vec{n}|-1}} w(x;\alpha-1,\beta_j-1|q) Q_{|\vec{n}|+1}(x),
\end{equation}
where $Q_{|\vec{n}|+1}$ is a monic polynomial of degree $|\vec{n}|+1$. We will show that this monic polynomial
$Q_{|\vec{n}|+1}$ satisfies the multiple orthogonality conditions (\ref{eq:2kind}) of
 $p_{\vec{n}+\vec{e}_j}(x;\alpha-1,\vec{\beta}-\vec{e}_j|q)$ and hence, since all $\beta_i-\beta_j \notin \mathbb{Z}$ whenever $i\neq j$, the unicity of the
multiple orthogonal polynomials implies that
$Q_{|\vec{n}|+1}(x) = p_{\vec{n}+\vec{e}_j}(x;\alpha-1,\vec{\beta}-\vec{e}_j|q)$.
Integration by parts gives
\begin{multline*}
  \frac{1-q^{\alpha+\beta_j+|\vec{n}|}}{(1-q)q^{\alpha+|\vec{n}|-1}}
  \int_0^1 x^k w(x;\alpha-1,\beta_j-1|q) Q_{|\vec{n}|+1}(x)\, d_qx   \\
 =   - q \int_0^1 w(x;\alpha,\beta_j|q) p_{\vec{n}}(x;\alpha,\vec{\beta}|q) \D_q x^k \, d_qx ,
\end{multline*}
so that
\[  \int_0^1 x^k w(x;\alpha-1,\beta_j-1|q) Q_{|\vec{n}|+1}(x)\, d_qx = 0, \qquad k=0,1,\ldots,n_j. \]
For the other components $\beta_i$ $(i \neq j)$ of $\vec{\beta}$ we have
\begin{eqnarray*}
  \lefteqn{ \frac{1-q^{\alpha+\beta_j+|\vec{n}|}}{(1-q)q^{\alpha+|\vec{n}|-1}}
  \int_0^1 x^k w(x;\alpha-1,\beta_i|q) Q_{|\vec{n}|+1}(x)\, d_qx} \quad & & \\
& = &    \frac{1-q^{\alpha+\beta_j+|\vec{n}|}}{(1-q)q^{\alpha+|\vec{n}|-1}}
 \int_0^1 x^k \frac{(q^{\beta_j}x;q)_\infty}{(q^{\beta_i+1}x;q)_\infty} w(x;\alpha-1,\beta_j-1|q) Q_{|\vec{n}|+1}(x)\, d_qx  \\
& = &   -\ q \int_0^1 w(x;\alpha,\beta_j|q) p_{\vec{n}}(x;\alpha,\vec{\beta}|q)
\D_q \left( x^k  \frac{(q^{\beta_j}x;q)_\infty}{(q^{\beta_i+1}x;q)_\infty}      \right)  \, d_qx ,
\end{eqnarray*}
and since $\beta_i-\beta_j \notin \mathbb{Z}$ we have
\[   \D_q \left( x^k \frac{(q^{\beta_j}x;q)_\infty}{(q^{\beta_i+1}x;q)_\infty} \right)
 = x^{k-1} \frac{(q^{\beta_j+1}x;q)_\infty}{(q^{\beta_i+1}x;q)_\infty} a_k(x) , \]
where each $a_k$ is a polynomial of degree exactly 1 and $a_0(0)=0$. Therefore
\[  \int_0^1 x^k w(x;\alpha-1,\beta_i|q) Q_{|\vec{n}|+1}(x)\, d_qx = 0, \qquad k=0,1,\ldots,n_i-1. \]
Hence all the orthogonality conditions for $p_{\vec{n}+\vec{e}_j}(x;\alpha-1,\vec{\beta}-\vec{e}_j|q)$ are
indeed satisfied. \qed

As a consequence we again find a \textit{Rodrigues formula}:

\begin{theorem}
The multiple little $q$-Jacobi polynomials of the second kind are given by
\begin{equation}  \label{eq:RodII}
 p_{\vec{n}}(x;\alpha,\vec{\beta}|q)  =
  \frac{C(\vec{n},\alpha,\vec{\beta})}{(qx;q)_\infty x^{\alpha}}
   \prod_{j=1}^r \left( (q^{\beta_j+1}x;q)_\infty \D_p^{n_j} \frac{1}{(q^{\beta_j+n_j+1}x;q)_\infty} \right)
    (qx;q)_\infty x^{\alpha+|\vec{n}|} ,
\end{equation}
where the product of the difference operators can be taken in any order and
\[  C(\vec{n},\alpha,\vec{\beta}) = (-1)^{|\vec{n}|}
    \frac{(1-q)^{|\vec{n}|} q^{(\alpha+|\vec{n}|-1)|\vec{n}|}}
    {\prod_{j=1}^r (q^{\alpha+\beta_j+|\vec{n}|+1};q)_{n_j}} . \]
\end{theorem}

\textbf{Proof:}
The proof can be given in a similar way as in the case of little $q$-Jacobi polynomials
of the first kind by repeated application of the raising operators. Alternatively one
can use induction on $r$. For $r=1$ the Rodrigues formula is the same as (\ref{eq:Rod}).
Suppose that the Rodrigues formula (\ref{eq:RodII}) holds for $r-1$. Observe that the multiple
orthogonal polynomials with multi-index $(n_1,\ldots,n_{r-1})$ for $r-1$ measures
$(\mu_1,\ldots,\mu_{r-1})$ coincide
with the multiple orthogonal polynomials with multi-index $(n_1,n_2,\ldots,n_{r-1},0)$
for $r$ measures $(\mu_1,\ldots,\mu_r)$ for any measure $\mu_r$. Use the Rodrigues
formula for $r-1$ for the polynomial
$p_{\vec{n}-n_r\vec{e}_r}(x;\alpha+n_r,\vec{\beta}+n_r\vec{e}_r|q)$ to find
\begin{multline*}
    w(x;\alpha+n_r,\beta_r+n_r|q) p_{\vec{n}-n_r\vec{e}_r}(x;\alpha+n_r,\vec{\beta}+n_r\vec{e}_r|q)
  =   C(\vec{n}-n_r\vec{e}_r,\alpha+n_r,\vec{\beta})
     \\  \times
       \frac{1}{(q^{\beta_r+n_r+1}x;q)_\infty}
   \prod_{j=1}^{r-1} \left( (q^{\beta_j+1}x;q)_\infty \D_p^{n_j} \frac{1}{(q^{\beta_j+n_j+1}x;q)_\infty} \right)
    (qx;q)_\infty x^{\alpha+|\vec{n}|}.
\end{multline*}
Now apply the raising operation (\ref{eq:rais2}) for $\beta_r$ to this
expression $n_r$ times to find the required expression.
\qed

In a similar way as for the first kind multiple little $q$-Jacobi polynomials we can
find an explicit formula with finite sums using the Rodrigues formula for little $q$-Jacobi
polynomials $r$ times. For $r=2$ this gives the following:

\begin{theorem}
The multiple little $q$-Jacobi polynomials of the second kind (for $r=2$) are explicitly
given by
\begin{multline}   \label{eq:explII}
p_{n,m}(x;\alpha,(\beta_1,\beta_2)|q) = \frac{q^{\alpha(n+m)+n^2+m^2+nm}(q^{-m-\alpha};q)_m(q^{-n-\alpha};q)_n
   (q^{\alpha+1};q)_{m+n}}{(q^{\alpha+\beta_1+n+m+1};q)_n(q^{\alpha+\beta_2+n+m+1};q)_m
   (q^{\alpha+1};q)_n(q^{\alpha+1};q)_m} \\
\times  \sum_{\ell=0}^n\sum_{k=0}^m \frac{(q^{-n};q)_\ell (q^{-m};q)_k(q^{\alpha+\beta_2+n+m+1};q)_k
      (q^{\alpha+\beta_1+n+1};q)_{k+\ell}}{(q^{\alpha+1};q)_{k+\ell} (q^{\alpha+\beta_1+n+1};q)_k} \\
  \times \frac{q^{k+\ell} x^{k+\ell}}{q^{nk} (q;q)_k(q;q)_\ell} .
\end{multline}
\end{theorem}

\textbf{Proof:}
The Rodrigues formula (\ref{eq:RodII}) for $r=2$ becomes
\begin{multline*}
p_{n,m}(x;\alpha,(\beta_1,\beta_2)|q) = \frac{(-1)^{n+m}(1-q)^{n+m} q^{(\alpha+n+m-1)(n+m)}}
   {(q^{\alpha+\beta_1+n+m+1};q)_n (q^{\alpha+\beta_2+n+m+1};q)_m} \\
\times  x^{-\alpha} \frac{(q^{\beta_1+1}x;q)_\infty}{(qx;q)_\infty}
  \D_p^n \frac{(q^{\beta_2+1}x;q)_\infty}{(q^{\beta_1+n+1}x;q)_\infty}
  \D_p^m \frac{(qx;q)_\infty}{(q^{\beta_2+m+1}x;q)_\infty} x^{\alpha+n+m} .
\end{multline*}
The Rodrigues formula (\ref{eq:Rod}) for little $q$-Jacobi polynomials gives
\begin{multline*}
  \D_p^m \frac{(qx;q)_\infty}{(q^{\beta_2+m+1}x;q)_\infty} x^{\alpha+n+m} \\
  = \frac{(-1)^m (q^{\alpha+\beta_2+n+m+1};q)_m}{(1-q)^m q^{\alpha m +m^2-m+nm}} x^{\alpha+n}
   \frac{(qx;q)_\infty}{(q^{\beta_2+1}x;q)_\infty} p_m(x;\alpha+n,\beta_2|q),
\end{multline*}
hence
\begin{multline*}
 p_{n,m}(x;\alpha,(\beta_1,\beta_2)|q) = \frac{(-1)^{n}(1-q)^{n} q^{\alpha n + n^2+nm-n}}
   {(q^{\alpha+\beta_1+n+m+1};q)_n } \\
\times  x^{-\alpha} \frac{(q^{\beta_1+1}x;q)_\infty}{(qx;q)_\infty}
  \D_p^n x^{\alpha+n} \frac{(qx;q)_\infty}{(q^{\beta_1+n+1}x;q)_\infty} p_m(x;\alpha+n,\beta_2|q).
\end{multline*}
Now use the explicit expression (\ref{eq:expl}) for the little $q$-Jacobi polynomials to find
\begin{multline*}
 p_{n,m}(x;\alpha,(\beta_1,\beta_2)|q) = \frac{(-1)^n(1-q)^n q^{\alpha (n+m) + n^2+m^2+2nm-n}(q^{-m-n-\alpha};q)_m}
   {(q^{\alpha+\beta_1+n+m+1};q)_n (q^{\alpha+\beta_2+n+m+1};q)_m } \\
\times  \frac{(q^{\beta_1+1}x;q)_\infty}{x^\alpha (qx;q)_\infty}
  \sum_{k=0}^m \frac{(q^{-m};q)_k(q^{\alpha+\beta_2+n+m+1};q)_k q^k}{(q^{\alpha+n+1};q)_k (q;q)_k}
 \D_p^n x^{\alpha+n+k} \frac{(qx;q)_\infty}{(q^{\beta_1+n+1}x;q)_\infty} .
\end{multline*}
Again we recognize a little $q$-Jacobi polynomial
\begin{multline*}
 \D_p^n x^{\alpha+n+k} \frac{(qx;q)_\infty}{(q^{\beta_1+n+1}x;q)_\infty} \\
 = \frac{(-1)^n (q^{\alpha+\beta_1+k+n+1};q)_n}{(1-q)^n q^{\alpha n + kn +n^2-n}}
  x^{\alpha+k} \frac{(qx;q)_\infty}{(q^{\beta_1+1}x;q)_\infty} p_n(x;\alpha+k,\beta_1|q),
\end{multline*}
and if we use the explicit expression (\ref{eq:expl}) for this little $q$-Jacobi polynomial, then
we find (\ref{eq:explII}) after some simplifications.
\qed

\section{Zeros}
The zeros of the multiple little $q$-Jacobi polynomials (first and
second kind) are all real, simple and in the interval $(0,1)$.
This is a consequence of the fact that $\mu_1,\ldots,\mu_r$ form
an AT-system \cite[first Corollary on p.~141]{NS}. For the usual
orthogonal polynomials with positive orthogonality measure $\mu$
we know that an interval $[c,d]$ for which the orthogonality
measure has no mass, i.e., $\mu([c,d])=0$, can have at most one
zero of each orthogonal polynomial $p_n$. In particular this means
that each orthogonal polynomial $p_n$ on the exponential lattice
$\{q^k,\ k=0,1,2,\ldots \}$ can have at most one zero between two
points $q^{k+1}$ and $q^k$ of the lattice. A similar result holds
for multiple orthogonal polynomials if we impose some conditions
on the measures $\mu_i$.

\begin{theorem}
Suppose $\mu_1,\ldots,\mu_r$ are positive measures on $[a,b]$ with infinitely many points in their support, which form an AT system,
i.e., $\mu_k$ is absolutely continuous with respect to $\mu_1$ for $2 \leq k \leq r$ with
\[  \frac{d\mu_k(x)}{d\mu_1(x)} = w_k(x), \]
and
\[  1,x,\ldots,x^{n_1-1},w_2(x),xw_2(x),\ldots,x^{n_2-1}w_2(x),\ldots,
w_r(x),xw_r(x),\ldots,    x^{n_r-1}w_r(x)
\]
are a Chebyshev system on $[a,b]$ for every multi-index $\vec{n}$. If $[c,d]$ is an interval such that $\mu_1([c,d])=0$, then each multiple orthogonal polynomial $p_{\vec{n}}$
has at most one zero in $[c,d]$.
\end{theorem}

\textbf{Proof:}
Suppose that $p_{\vec{n}}$ is a multiple orthogonal polynomial with two zeros
$x_1$ and $x_2$ in $[c,d]$. We can then write it as
$p_{\vec{n}}(x)=(x-x_1)(x-x_2)q_{|\vec{n}|-2}(x)$, where $q_{|\vec{n}|-2}$ is a polynomial
of degree $|\vec{n}|-2$.
Consider a function $A(x) = \sum_{j=1}^r A_j(x)w_j(x)$, where $w_1=1$ and each $A_j$ is a polynomial
of degree $m_j-1\leq n_j-1$, with $|\vec{m}|=|\vec{n}|-1$. Since we are dealing with a Chebyshev system, there is a unique function $A$ satisfying the interpolation conditions
\[   A(y) = \begin{cases} 0 & \quad \textrm{if $y$ is a zero of $q_{|\vec{n}|-2}$}, \\
                  1 & \quad \textrm{if $y = x_1$}.
        \end{cases}  \]
Furthermore $A$ has $|\vec{n}|-2$ zeros in $[a,b]$ and these are the only
sign changes on $[a,b]$. Hence
\[  \int_a^b p_{\vec{n}}(x) A(x) \, d\mu_1(x) =
    \int_{[a,b]\setminus [c,d]} (x-x_1)(x-x_2) q_{|\vec{n}|-2}(x) A(x)\, d\mu_1(x) \neq 0, \]
since the integrand does not change sign on $[a,b] \setminus [c,d]$.
On the other hand
\[  \int_a^b p_{\vec{n}}(x) A(x) \, d\mu_1(x) =
    \sum_{j=1}^r \int_a^b p_{\vec{n}}(x)A_j(x) \,d\mu_j(x) = 0, \]
since every term in the sum vanishes because of the orthogonality conditions. This contradiction
implies that $p_{\vec{n}}$ can't have two zeros in $[c,d]$.
\qed

In particular this theorem tells us that the zeros of the multiple little $q$-Jacobi polynomials
are always separated by the points $q^k$ and that between two points $q^{k+1}$ and $q^k$ there
can be at most one zero of a multiple little $q$-Jacobi polynomials. Note that the points
$q^k$ have one accumulation point at 0, hence as a consequence the zeros of the
multiple little $q$-Jacobi polynomials (first and second kind) accumulate at the origin.

\section{Asymptotic behavior}
The asymptotic behavior of little $q$-Jacobi polynomials was given by Ismail and Wilson \cite{IW} and an asymptotic expansion was given by Ismail \cite{I}. In this section we
give the asymptotic behavior of the multiple little $q$-Jacobi polynomials which extends
the result of Ismail and Wilson.

\begin{theorem}
For the multiple little $q$-Jacobi polynomials of the first kind we have
\begin{equation} \label{eq:1kindas}
  \lim_{n,m \to \infty}  x^{n+m} p_{n,m}(1/x;(\alpha_1,\alpha_2),\beta|q) = (x;q)_\infty.
\end{equation}
The order in which the limits for $n$ and $m$ are taken is irrelevant.
\end{theorem}
\textbf{Proof:}
If we use (\ref{eq:explI}) and reverse the order of summation (i.e., change variables $m-k=j$ and $n-\ell=i$),
then
\begin{multline*}
  x^{n+m}p_{n,m}(1/x;(\alpha_1,\alpha_2),\beta|q) = \frac{q^{nm +m^2+n^2+\alpha_1n+\alpha_2m}
   (q^{-\alpha_1-n};q)_n(q^{-\alpha_2-m};q)_m}{(q^{\alpha_1+\beta+n+m+1};q)_n(q^{\alpha_2+\beta+n+m+1};q)_m} \\
\times  \sum_{i=0}^n \sum_{j=0}^m \frac{(q^{-n};q)_{n-i} (q^{-m};q)_{m-j} (q^{\alpha_2+\beta+m+n+1};q)_{m-j}
(q^{\alpha_1+\beta+n+1};q)_{m+n-i-j} (q^{\alpha_1+n+1};q)_{m-j}}{(q^{\alpha_2+1};q)_{m-j}(q^{\alpha_1+1};q)_{m+n-i-j}
(q^{\alpha_1+\beta+n+1};q)_{m-j}} \\
 \times \frac{q^{m+n-i-j} x^{i+j}}{q^{(m-j)n}(q;q)_{m-j} (q;q)_{n-i}} .
\end{multline*}
Now observe that
\begin{eqnarray*}
  (q^{-m};q)_{m-j} & = & (-1)^{m-j} q^{-\frac{m(m+1)}2+\frac{j(j+1)}2} \frac{(q;q)_m}{(q;q)_j}, \\
  (q^{-m-\alpha};q)_m & = & (-1)^m q^{-m(m+1)/2} q^{-m\alpha} (q^{\alpha+1};q)_m, \\
  (q^{c+n};q)_m & = & \frac{(q^c;q)_{n+m}}{(q^c;q)_n},
\end{eqnarray*}
therefore we find
\begin{multline*}
  x^{n+m}p_{n,m}(1/x;(\alpha_1,\alpha_2),\beta|q) = \frac{(q^{\alpha_2+1};q)_m
   (q^{\alpha_1+\beta+1};q)_{n+m} (q;q)_m(q;q)_n}
   {(q^{\alpha_1+\beta+1};q)_{2n+m} (q^{\alpha_2+\beta+1};q)_{n+2m}} \\
\times  \sum_{i=0}^n \sum_{j=0}^m  \frac{(q^{\alpha_2+\beta+1};q)_{n+2m-j} (q^{\alpha_1+\beta+1};q)_{2n+m-i-j}
  (q^{\alpha_1+1};q)_{n+m-j}}
  { (q^{\alpha_2+1};q)_{m-j} (q^{\alpha_1+1};q)_{m+n-i-j}   (q^{\alpha_1+\beta+1};q)_{m+n-j}(q;q)_{n-i}(q;q)_{m-j}} \\
 \times (-1)^{i+j} q^{\binom{i}{2}+\binom{j}{2}} \frac{q^{nj} x^{i+j}}{(q;q)_i(q;q)_j}.
\end{multline*}
If we use Lebesgue's dominated convergence theorem, then we take $n,m \to \infty$
in each term of the sum. The factor $q^{nj}$ tends to zero whenever $j > 0$, hence
the only contributions come from $j=0$, and we find
\[  \lim_{n,m\to \infty} x^{n+m}p_{n,m}(1/x;(\alpha_1,\alpha_2),\beta|q)
  = \sum_{i=0}^\infty q^{\binom{i}{2}} \frac{(-x)^i}{(q;q)_i}. \]
The right hand side is the $q$-exponential function
\[  E_q(-x) = (x,q)_\infty, \]
\cite[(II.2) in Appendix II]{GR}, which gives the required result.
\qed

\begin{theorem}
For the multiple little $q$-Jacobi polynomials of the second kind we have
\begin{equation} \label{eq:2kindas}
  \lim_{n,m \to \infty}  x^{n+m} p_{n,m}(1/x;\alpha,(\beta_1,\beta_2)|q) = (x;q)_\infty.
\end{equation}
The order in which the limits for $n$ and $m$ are taken is irrelevant.
\end{theorem}
\textbf{Proof:}
The proof is similar to the case of the first kind multiple little $q$-Jacobi polynomials,
except that now we use the expression (\ref{eq:explII}).
\qed

As a consequence (using Hurwitz' theorem) we see that every zero of $(1/x;q)_\infty$, i.e., each number $q^k,\ k=0,1,2,\ldots$, is an accumulation point of zeros of the multiple little $q$-Jacobi
polynomial $p_{n,m}$ of the first and of the second kind.

\end{document}